\documentclass{amsart}
\usepackage{amsmath,amsfonts,amssymb,amscd,delarray,verbatim,epigraph}

\DeclareMathSymbol{\twoheadrightarrow}  {\mathrel}{AMSa}{"10}

\def\Q{{\mathbb Q}}
\def\Z{{\mathbb Z}}

\def\Gal{\mathrm{Gal}}

\def\Br{\mathrm{Br}}
\def\TORS{\mathrm{TORS}}
\def\DIV{\mathrm{DIV}}
\def\et{\mathrm{et}}
\def\Gm{\mathbb{G}_m}
\def\NS{\mathrm{NS}}
\def\Pic{\mathrm{Pic}}
                                                                            
\def\Hom{\mathrm{Hom}}

\def\s{{\mathfrak S}}

\theoremstyle{definition}

\title[The Brauer Group of a Surface over a Finite Field]
{The Brauer Group of a Surface over a Finite Field}
\author[Yu.\ G.\ Zarhin]{Yu.\ G.\ Zarhin}

\begin{document}



\maketitle

Let $k$ be a finite field of characteristic $p$ that consists of $q$ elements, $ \bar{k}$ its algebraic closure, $\Gamma=\Gal(\bar{k}/k)$  its absolute Galois group,
$$\sigma:  \bar{k} \to  \bar{k}, \ a \mapsto a^q$$
the Frobenius automorphism, which is a canonical generator of procyclic $\Gamma$. Let $X$ be an absolutely irreducible smooth projective variety over $k$, $\Br^{\prime}(X)=H^2(X_{\et},\Gm)$ its cohomological Brauer group; if $X$ is a surface then $\Br^{\prime}(X)$ coincides with the Brauer group $\Br(X)$, which is defined as the set of similarily classes of Azumaya algebras over $X$ \cite{2,7}. The group $\Br^{\prime}(X)$ is periodic commutative; for each positive integer $n$ the kernel $_{n}\Br^{\prime}(X)$ of multiplication by $n$ in $\Br^{\prime}(X)$ is finite \cite{5,7}.
The Artin-Tate conjecture  asserts that $\Br^{\prime}(X)$ is always finite; this conjecture is closely related to the Tate conjecture that deals with the order of pole of the zeta function of $X$ at $1$  \cite{2,5,9,15,1}. Let us put $\bar{X}=X\otimes \bar{k}$. If $X$ is a surface then the finiteness of
 $\Br(X)=\Br^{\prime}(X)$ is known in  each of the following cases \cite{5,7,6,11}.
 
 \begin{enumerate}
  \item
  $X$ is an abelian surface;
   \item
   $\bar{X}$ is a rational surface;
    \item
    $\bar{X}$ is birational to a product of two curves;
     \item
     $X$ is either a Kummer surface or a K3 surface with a pencil of elliptic curves
     \item
     $X$ is either ordinary or a supersingular K3 surface.
     \end{enumerate}

     If $X$ is a surface then there is a nondegenerate skewsymmetric bilinear pairing \cite{5,9}
     $$\Br(X)_{\DIV} \times \Br(X)_{\DIV}  \to \Q/\Z.$$
     Hereafter if $B$ is any commutative group then we write
     $B_{\DIV}=B/\DIV  B$ where 
     $$\DIV  B=\mathrm{Im}[\Hom(\Q,B) \to \Hom(\Z,B)=B]$$
     is the subgroup of infinitely divisible elements of $B$. (By definition,
     $\DIV  B \subset \bigcap_n nB$; if for each $n$ the kernel $_{n}B$ of multiplication by $n$ is finite the the inclusion becomes equality.)
     
     According to Tate \cite{5}, a question whether this pairing is alternating is an ``interesting cohomological problem''. (Recall that the skewsymmetry  means that the value of the pairing multiplies by $-1$ if one permutes the arguments while the alternation means the vanishing of the value of the pairing if the arguments do coincide.) Of course, it suffices to check the alternation on the $2$-primary component  $\Br(X)_{\DIV}(2)=\Br(X)(2)_{\DIV}$ of $\Br(X)_{\DIV}$.
     Hereafter if $B$ is a commutative group and $l$ is a prime then we write $B(l)$ for its subgroup of all elements, whose order is a power of $l$. We have
     $$\oplus B(l)=\TORS(B)=\oplus \TORS(B)(l)$$
       where $\TORS(B)$ is the subgroup of all elements of finite order in $B$ and the summation is taken over all primes $l$. If $B$ is a periodic group then $B=\TORS(B)$ and
        $$\DIV B=\oplus \DIV(B(l))=\oplus \DIV(B)(l),$$
     $$B=\oplus B(l), \  B_{\DIV}=\oplus  B(l)_{\DIV}=\oplus  B_{\DIV}(l)$$
   and the summation is taken over all primes $l$. If the group $_{l}B$ is finite then the group $B_{\DIV}(l)$ is also finite. This implies that all the groups $\Br^{\prime}(X)_{\DIV}(l)$ are finite.  If $X$ is a surface then  the order of  $\Br(X)_{\DIV}(l)$ is a square for all odd $l$. The alternation of the pairing above would imply that the order of the $2$-primary component $\Br(X)_{\DIV}(2)$ is  also a square.
   
   Of course,  if   $\Br(X)$ is finite then  $\Br(X)= \Br(X)_{\DIV}$. It is known \cite{1} that $\Br(X)_{\DIV}$ is finite if $X$ is a K3 surface
   
   The aim of this note is to prove thar the order of $\Br(X)_{\DIV}(2)$ is a square under the following additional assumptions on the surface $X$.
   
   \begin{enumerate}
   \item
   The surface $\bar{X}=X\otimes \bar{k}$ lifts to characteristic zero;
   \item
   $p \ne 2$;
   \item
   There is no $2$-torsion in the N\'eron-Severi group $\NS(\bar{X})$ of $\bar{X}$.
\end{enumerate}

The proof is based on the construction of a nondegenerate skewsymmetric bilinear pairing
$$\Br(X)_{\DIV}(l)\times \Br(X)_{\DIV}(l) \to \Q/\Z$$
under assumptions that $l \ne p$ and there is no $l$-torsion in $\NS(\bar{X})$.  When $p \ne 2$ we prove that if $X$ is liftable to characteristic zero then this pairing is alternating for all $l$, including the case $l=2$. It would be interesting to compare this pairing with the pairing constructed in \cite{5}.

Notice that the present paper is a nataral complement to our paper \cite{1} despite of slightly different notation.

{\bf A construction of the pairing}.  In the next lemma and after it,  if $A$ is a commutative group provided with a structure of a $\Gamma$-module then we write $A^{\Gamma}$ and $A_{\Gamma}$ for its subgroup of $\Gamma$-invariants and the quotient group of $\Gamma$-coinvariants respectively. i.e.
$$A^{\Gamma}=\{a \in A\mid \sigma a=a\}, \ A_{\Gamma}=A/(1-\sigma)A.$$
All the cohomology groups of algebraic varieties in this paper are taken with respect to \'etale topology.

{\bf Lemma 1}. {\sl Let $X$ be a surface, $l \ne p$ and $\NS(\bar{X})(l)=0$.

Then $\TORS(H^2(\bar{X},\Z_l(1)))=0 \ $,  $\ H^2(\bar{X},\Z_l(1))$ is a free $\Z_l$-module of finite rank and there is a natural isomorphism}
$$\Br(X)_{\DIV}(l)=\TORS[ H^2(\bar{X},\Z_l(1))_{\Gamma}].$$

\begin{proof}
Recall \cite{5,7} that there is a natural embedding
$$\NS(\bar{X})\otimes \Z_l \to H^2(\bar{X},\Z_l(1))$$
that induces an isomorphism of the corresponding torsion subgroups \cite[III, Sect. 8.2]{7}.
So, if finitely generated group $\NS(\bar{X})$ has no $l$-torsion then finitely generated $\Z_l$-module $H^2(\bar{X},\Z_l(1))$ has no torsion as well and therefore is free. So,
 $$\TORS(H^2(\bar{X},\Z_l(1)))=0.$$
 According to proposition 1.4.2 of \cite{1}, there is an exact sequence
 $$0 \to \TORS[H^2(\bar{X},\Z_l(1))_{\Gamma}]\to \Br(X)_{\DIV}(l) \to \TORS[H^3(\bar{X},\Z_l(1))]^{\Gamma}.$$
 In order to finish the proof, it suffices to check that $\TORS[H^3(\bar{X},\Z_l(1))]=0$.
 The Poincar\'e duality for $\bar{X}$ implies that finite commutative $l$-groups $\TORS(H^2(\bar{X},\Z_l(1)))$ and $\TORS(H^3(\bar{X},\Z_l(1)))$ are mutually dual   \cite[III]{7}. Therefore  $\TORS(H^3(\bar{X},\Z_l(1)))$ is also zero.
\end{proof}

    {\bf Lemma 2}. {\sl Let $X$ be a surface, $l \ne p$ and $\NS(\bar{X})=0$. Then 
 $H^2(\bar{X},\Z_l(1))$  is a free $\Z_l$-module of finite rank  and the intersection pairing
 $$<,>:  H^2(\bar{X},\Z_l(1)) \times H^2(\bar{X},\Z_l(1)) \to Z_l$$
 is perfect, i.e., the corresponding homomorphism of  free $\Z_l$-modules
 $$H^2(\bar{X},\Z_l(1)) \to \Hom_{\Z_l}(H^2(\bar{X},\Z_l(1)),\Z_l)$$
 is an isomorphism. The form $<,>$ is $\Gamma$-invariant, i.e.}
 $$,\sigma x,\sigma y>=<x,y> \ \forall x,y \in H^2(\bar{X},\Z_l(1)).$$

 \begin{proof}
 It is well known \cite{5} that $<,>$ is $\Gamma$-invariant. By Lemma 1, 
 $H^2(\bar{X},\Z_l(1))$ is a free $\Z_l$-module. It follows from the proof of Lemma 1
 that $H^3(\bar{X},\Z_l(1))$ is rorsion-free. Hence the natural embedding
 $$H^2(\bar{X},\Z_l(1))\otimes\Z/l\Z \to H^2(\bar{X},\mu_l)$$
 is an isomorphism. Here $\mu_l$ is the sheaf of $l$th roots of unity on $\bar{X}_{\et}$ \cite{2}.
 Now the perfectness of $<,>$ follows readily from Poincar\'e duality for  $\mu_l$ \cite{2} and Nakayama's Lemma. (Compare with arguments in the last section of Schneider's paper \cite{14}.)
 \end{proof}

{\bf Lemma 3}. {\sl Let $H$ be a free $\Z_l$-module of finite rank provided with a continuous action of $\Gamma$,  
$$(,): H \times H \to \Z_{\ell}$$
a perfect $\Gamma$-invariant symmetric pairing, i.e.,
$$(\sigma x,\sigma y)=(x,y) \ \forall x,y \in H$$
and the homomorpism of free $\Z_l$-modules 
$$H \to \Hom_{\Z_l}(H,\Z_l)$$
induces by $(,)$ is an isomorphism.

Then there exists a natural nondegenerate skewsymmetric bilinear pairing
$$(,)_B: \TORS(H_{\Gamma}) \times  \TORS(H_{\Gamma}) \to \Q_{l}/\Z_{l}\subset \Q/\Z.$$
In particular, if $l \ne 2$ then the order of finite abelian group $\TORS(H_{\Gamma})$ is a full square.

Let $l=2$ and $w \in H$ a characteristic element of $(,)$, i.e.
$$(x,x)-(x,w)\in 2 \Z_2 \ \forall x \in H.$$

If $w$ is $\Gamma$-invariant, i.e., $w \in H^{\Gamma}$, then $(,)_B$ is alternating; in particular,
the order of $\TORS(H_{\Gamma})$ is a full square.

If there is an odd positive integer $n$ such that $\sigma^n w=w$ then $(,)_B$ is also alternating and
the order of $\TORS(H_{\Gamma})$ is a full square.}

\begin{proof}
We choose as $(,)_B$ the pairing constucted in subsections 3.3, 3.4 of \cite{1}.
(In these subsections $\TORS(H_{\Gamma})$ is denoted by $B_0(H)$ and $H^{\Gamma}$ by $H^{\sigma}$.) It is proven (ibid) that $(,)_B$ is nondegenerate skewsymmetric.  There is also the following alternation criterion  of  $(,)_B$  \cite[Sect. 3.4.1]{1}.

{\sl Let $l=2$. Denote by $H^0$ the orthogonal complement to $H^{\Gamma}$ with respect to $(,)$, i.e.,
$$H^0:=\{x\in H\mid (x,y)=0 \ \forall y \in H^{\Gamma}\}.$$
The pairing $(,)_B$ is alternating if and only if the restriction of $(,)$ to $H^0$ is even, i.e.,}
$$(x,x)\in 2\Z_2 \ \forall x \in H^0.$$

Let us assume that the characteristic element $w \in H^{\Gamma}$. Then $(x,w)=0$ for all $x \in H^0$ and therefore, by definition of a characteristic element,
$$(x,x)\in 2\Z_2 \ \forall x \in H^0.$$
So, $(,)_B$ satisfies the alternation criterion.

If $\sigma^n w=w$ for some odd positive integer $n$ then one may easily check that
$$w^{\prime}=w+\sigma w+ \dots + \sigma^{n-1}w$$
is a $\Gamma$-invariant characteristic element and apply  the same alternation criterion.
\end{proof}

{\bf Remark}.
For reader's convenience we reproduce an explicit construction of $(,)_B$. In order to do that, let us put
$H^{\prime}=H\otimes_{\Z_l}\otimes \Q_{\ell}$ and extend the form $(,)$ by $\Q_l$-linearity to the symmetric $\Q_l$-bilinear pairing
$$H^{\prime}\times H^{\prime} \to  \Q_{\ell},$$
which we continue to denote $(,)$.  Then for each
$$a,b \in \TORS(H_{\Gamma})\subset H/(1-\sigma)H$$
we have
$$(a,b)_B:=(u,v)\bmod \Z_l \in \Q_l/\Z_l$$
where $u \in H, v \in H^{\prime}$ satisfy
$$a=u\bmod (1-\sigma)H, \ v-\sigma v\in H, \ b=(v-\sigma v) \bmod (1-\sigma)H.$$

\vskip .2cm

Lemmas 1 and 2 combined with Lemma 3 applied to $H=H^2(\bar{X},\Z_l(1))$ and $(,)=<,>$
give us the following assertion.

\vskip .2cm

{\bf Theorem 1}.  {\sl Let $X$ be a surface, $l \ne p$ and $\NS(\bar{X})(l)=0$. 

Then $\TORS(H^2(\bar{X},\Z_l(1)))=0$ and there exists a natural nondegenerate skewsymmetric bilinear pairing
$$<,>_B: \Br(X)_{\DIV}(l) \times  \Br(X)_{\DIV}(l)  \to \Q_l/\Z_l\subset \Q/\Z.$$
In particular, if $l \ne 2$ then the order of finite abelian group $\Br(X)_{\DIV}(l)$ is a full square. 

If $l=2$ and there exists $w \in H^2(\bar{X},\Z_l(1))^{\Gamma}$ such that
$$(x,x)-(x,w) \in 2 \Z_2 \ \forall x \in H^2(\bar{X},\Z_2(1))$$
then $<,>_B$ is alternating. In particular,  the order of finite abelian group $\Br(X)_{\DIV}(2)$ is a full square.}

\vskip .2cm

{\bf Example}. Let us assume that either $X$ is a K3 surface or an abelian surface, or $\bar{X}$ is biregular over $\bar{k}$ to a product of two curves. Then $\NS(\bar{X})(l)=0$ for all $l \ne p$.
If $p \ne 2$ then for $l=2$ the intersection pairing $<,>$ is even and one may put $w=0$. For abelian and K3 surfaces the evenness and the absense of torsion follows from the existence of a lifting in characteristic 0 \cite{9,10,12} and for products of curves from the 
K{u}nneth formula. So, under these assumptions on $X$ and $p \ne 2$ the natural nondegenerate bilinear pairing
$$<,>_B: \Br(X)_{\DIV}(l) \times  \Br(X)_{\DIV}(l)  \to \Q_l/\Z_l\subset \Q/\Z$$
is alternating for all primes $l \ne p$, including $l=2$. In particular, for all $l \ne p$ the order of 
$ \Br(X)_{\DIV}(l)$ is a full square. In the case $l =p$ (and $p\ne 2$) Milne \cite{8} constructed
a nondegenerate bilinear skewsymmetric bilinear pairing
$$<,>_B: \Br(X)_{\DIV}(p) \times  \Br(X)_{\DIV}(p)  \to \Q_p/\Z_p\subset \Q/\Z,$$
which is automatically alternating, since $p\ne 2$.

Recall that $\Br(X)$ is finite if $X$ is an abelian variety or $\bar{X}$ is biregular to a product of two curves. If $X$ is a K3 surface then $\Br(X)_{\DIV}$ is a finite group. Combining all these assertions, we obtain the following statement.

\vskip .2cm

{\bf Theorem 2}. {\sl
Suppose that $p \ne 2$.
\begin{itemize}
\item[(a)] If $X$ is an abelian surface or $\bar{X}$ is biregular to a product of two curves then there exists a natural nondegenerate alternating pairing
$$\Br(X) \times \Br(X) \to \Q/\Z.$$
In particular, the order of finite abelian group $\Br(X)$ is a full square.
\item[(b)] 
If $X$ is a K3 surface then there exists a natural nondegenerate alternating pairing
$$\Br(X)_{\DIV} \times \Br(X)_{\DIV} \to \Q/\Z.$$
In particular, the order of finite abelian group $\Br(X)_{\DIV}$ is a full square. 
\end{itemize}}
\vskip .2cm

Let us assume that $\bar{X}$ lifts to characteristic $0$, i.e., there exist a  commutative complete discrete valuation ring $O$ with fraction field $K$ of characteristic zero and residue field $\bar{k}$ and a smooth projective morphism $f: \mathcal{Y} \to S=\mathrm{Spec }\ O$, whose closed fiber coincides with $\bar{X}$.  If $n$ is any positive integer not divisible by $p$ then there is an exact sequence
$$0 \to \mu_n \to \Gm \overset{n}{\to}  \Gm \to 0$$
of sheaves on $\mathcal{Y}_{\et}$,
which induces the following exact sequence 
$$R^{1}f_{*}\Gm \overset{n}{\to} R^{1}f_{*}\Gm \to R^{2}f_{*}\mu_n$$ of sheaves on $S_{\et}$,
which is part of the long exact sequence of higher direct images with respect to $f$.
Taking the corresponding groups of global sections, we get a canonical map
$$\delta: H^0(S, R^{1}f_{*}\Gm)\to  H^0(S, R^{2}f_{*}\mu_n).$$
Proper and smooth base change theorems \cite{2} imply that $R^{2}f_{*}\mu_n$ is a constant sheaf on $S_{\et}$, whose fiber over $\bar{k}$ coincides with $H^2(\bar{X}, \mu_n)$ while the fiber over $\bar{K}$ coincides with  $H^2(\bar{\mathcal{Y}}, \mu_n)$ where $\bar{K}$ is an algebraic closure of $K$ and $\bar{\mathcal{Y}}= \mathcal{Y}\otimes \bar{K}$ is the generic geometric fiber of $f$. In this, there are natural isomorphisms
$$H^2(\bar{\mathcal{Y}}, \mu_n)=H^0(S, R^{2}f_{*}\Gm) =H^2(\bar{X}, \mu_n).$$
Recall \cite{13} that there is a natural homomorphism from the Picard group
$$\s: \Pic(\mathcal{Y})=H^1(\mathcal{Y},\Gm) \to H^0(S, R^{1}f_{*}\Gm)$$
and also  homomorphisms of Picard groups \cite{2}
$$c_1: \Pic(\bar{X})=H^1(\bar{X},\Gm) \to H^2(\bar{X}, \mu_n),$$
$$c_1: \Pic(\bar{\mathcal{Y}})=H^1((\bar{\mathcal{Y}},\Gm) \to H^2((\bar{\mathcal{Y}}, \mu_n)$$
that arise from exact sequences of sheaves
$$0 \to \mu_n \to \Gm \overset{n}{\to} \Gm \to 0.$$
In this, the images of the canonical classes
$$w_{\bar{X},n}=c_1(\Omega^2_{\bar{X}}) \in H^2(\bar{X}, \mu_n), \
w_{\bar{\mathcal{Y}},n}=c_1(\Omega^2_{\bar{\mathcal{Y}}}) \in H^2(\bar{\mathcal{Y}},
 \mu_n)$$
 go to each other under the natural isomorphism
 $$H^2(\bar{\mathcal{Y}}, \mu_n)=H^0(S, R^{2}f_{*}\Gm) =H^2(\bar{X}, \mu_n).$$
 In order to prove it, it suffices to consider the homomorphism
 $$c_1:\Pic(\mathcal{Y}) =H^1(\mathcal{Y},\Gm) \overset{\s}{\to}  H^0(S, R^{1}f_{*}\Gm) \overset{\delta}{\to}  H^0(S, R^{2}f_{*}\mu_n)$$
 and notice that the images of $w_{\bar{X},n}$ and $w_{\bar{\mathcal{Y}},n}$ coincide with the image $c_1(\Omega^2_{\mathcal{Y}/S})$ of the relative canonical class $\Omega^2_{\mathcal{Y}/S}$, which follows immediately from the commutativeness of the  diagrams.
 
 If we take as $n$  various powers $l^i$ of $l$  and take  projective limits with respect to $i$, then we obtain the natural isomorphism of $\Z_{\ell}$-modules
 $$H^2(\bar{X},\Z_l(1)) =H^2(\bar{\mathcal{Y}},\Z_l(1)),$$
 which repects the intersection pairings. Elements $w_{\bar{X},n}$ and $w_{\bar{\mathcal{Y}},n}$  are compatible for various $n$ and are glueing together to the elements
 $$w_{\bar{X}} \in H^2(\bar{X},\Z_l(1)),  \  w_{\bar{\mathcal{Y}}}\in H^2(\bar{\mathcal{Y}},\Z_l(1)),$$
 which go to each other under the natural isomorphism above;  their respective images in
 $$ H^2(\bar{X},\Z_l(1))/\ell^i  H^2(\bar{X},\Z_l(1))\subset H^2(\bar{X},\mu_{l^i})$$
 and
 $$H^2(\bar{\mathcal{Y}},\Z_l(1))/\ell^i H^2(\bar{\mathcal{Y}},\Z_l(1))\subset 
 H^2(\bar{\mathcal{Y}},\mu_{\ell^i})$$
 coincide with $w_{\bar{X},\ell^i}$ and $w_{\bar{\mathcal{Y}},\ell^i}$ respectively.
 Since $\bar{X}$ is obtained from $k$-surface $X$ by extensions of scalars, all $w_{\bar{X},n}$  are $\Gamma$-invariant. This implies that
 $$w_{\bar{X}} \in H^2(\bar{X},\Z_l(1))^{\Gamma}.$$
 
 Now let us assume that $p>2, \ l=2$ and $\NS(\bar{X})(2)=0$. Then $H^2(\bar{X},\Z_2(1))$  and $H^2(\bar{\mathcal{Y}},\Z_2(1))$ are free $\Z_2$-modules. Wu's  theorem that relates Stiefel-Whitney classes and Steenrod squares \cite{3} and comparison theorems for clasical and \'etale cohomology \cite{2} imply that $w_{\bar{\mathcal{Y}}}$ is a characteristic element of the intersection pairing on $H^2(\bar{\mathcal{Y}},\Z_l(1))$. This implies that $w_{\bar{X}}$ is a characteristic element of the intersection pairing on $H^2(\bar{X},\Z_l(1))$. The $\Gamma$-invariance of  $w_{\bar{X}}$ combined with Theorem 1 imply the following result.
 
 \vskip .3cm
 
 {\bf Theorem 3}. {\sl Let $X$ be a surface, $p \ne 2$ and $\NS(\bar{X})(2)=0$.
 Let us assume that $\bar{X}$ lifts to characteristic $0$. Then there exists a natural nondegenerate alternating pairing
 $$\Br(X)_{\DIV}(2) \times \Br(X)_{\DIV}(2)  \to \Q_2/\Z_2\subset \Q/\Z.$$
 In particular, the order of finite abelian group $\Br(X)_{\DIV}(2)$ is a full square.} 
 
 \vskip .3cm
 
 {\bf Remark}. A question whether $w_{\bar{X}}$ is a characteristic element without assuming the existence of a lifting to characteristic $0$ was raised in \cite{4}.
 
  \vskip .3cm
 
 {\bf Remark}.  If $X$ lifts to characteristic $0$ then it is easy to see that the element  $w_{\bar{\mathcal{Y}}}$ is $\Gamma$-invariant and  characteristic one, which somehow simplifies the proof of Theorem 3.

\end{document}